\documentstyle[12pt,amssym]{article}

\newtheorem{theo}{Theorem}
\newtheorem{lemm}[theo]{Lemma}

\newcommand{\myF}{{\cal F}}

\newcommand{\myM}{{\cal M}}

\newcommand{\bR}{{\Bbb R}}

\newcommand{\trace}{{\rm Tr}}

\newcommand{\levy}{L\'{e}vy\,}

\begin{document}

\begin{center}
{\large\bf A three-series theorem on Lie groups}

Ming Liao\footnote{Department of Mathematics, Auburn University, Auburn, AL 36849, USA. Email: liaomin@auburn.edu}

\end{center}

\begin{quote}
{\bf Summary} \ We obtain a necessary and sufficient condition for the convergence of independent products on Lie groups, as a natural extension
of Kolmogorov's three-series theorem. Application to independent random matrices is discussed.



\end{quote}

\section{Introduction and main results} \label{sec1}

Let $x_n$ be a sequence of independent real-valued random variables.
Fix any constant $r>0$. Kolmogorov's three-series theorem (see
for example \cite[Theorem~22.8]{billingsley}) states that the series $\sum_{n=1}^\infty x_n$ converges
almost surely if and only if the following three conditions hold.
\vspace{2ex}

\noindent (K1) \ $\sum_{n=1}^\infty P(\vert x_n\vert>r)\ <\ \infty$;

\noindent (K2) \ $\sum_{n=1}^\infty E(x_n1_{[\vert x_n\vert\leq r]})$ converges, where $1_A$ is the indicator of a set $A$; and

\noindent (K3) \ $\sum_{n=1}^\infty E[(x_n1_{[\vert x_n\vert\leq r]}-b_n)^2]\ <\ \infty$, where $b_n=E(x_n1_{[\vert x_n\vert\leq r]})$.
\vspace{2ex}

Extensions of the three-series theorem to more general spaces have been explored in literature. In particular, Maksimov \cite{maksimov}
obtained a one-sided extension of the three-series theorem to Lie groups, providing a set of sufficient conditions for the convergence of
products of independent random variables in a Lie group, with some partial result toward the more difficult necessity part.

The purpose of this paper is to present a complete extension of the three-series theorem to a general Lie group. Our result
is a simpler form of a conjecture proposed in \cite{maksimov}, and is in more close analogy with the classical result.
We not only establish the more difficult necessity part, the proof
of sufficiency is also much shorter than \cite{maksimov}. The result will be applied to study the convergence of products
of independent random matrices.

Let $G$ be a Lie group of dimension $d$ with identity element $e$. There are a relatively compact neighborhood $U$ of $e$
and a smooth function $\phi=(\phi_1,\phi_2,\ldots,\phi_d)$: $U\to\bR^d$ which maps $U$ diffeomorphically onto
a convex neighborhood $\phi(U)$ of the origin $0$ in $\bR^d$, with $\phi(e)=0$. The $U$ is not assumed to be open
and $\phi$ is assumed extendable to be a smooth function on an open set containing the closure $\overline{U}$ of $U$. In the rest of the paper, $U$ and $\phi$ are fixed, but they may be chosen arbitrarily as long as the above properties are satisfied.

Let $x$ be a random variable in $G$. Its $U$-truncated mean $b$ is defined by
\begin{equation}
\phi(b) = E[\phi(x)1_{[x\in U]}]. \label{a}
\end{equation}
  Note that because $\phi(U)$ is convex, $E[\phi(x)1_{[x\in U]}]\in\phi(U)$ and $b=\phi^{-1}\{E[\phi(x)1_{[x\in U]}]\}$.

\begin{theo} \label{th1}
Let $x_n$ be a sequence of independent $G$-valued random variables with $U$-truncated means $b_n$. Then $\hat{x}_n=
x_1x_2\cdots x_n$ converges almost surely in $G$ as $n\to\infty$ if and only if the following three conditions hold.
\vspace{2ex}

\noindent (G1) \ $\sum_{n=1}^\infty P(x_n\in U^c)\ <\ \infty$, where $U^c$ is the complement of $U$ in $G$;

\noindent (G2) \ $\hat{b}_n=b_1b_2\cdots b_n$ converges in $G$ as $n\to\infty$; and

\noindent (G3) \ $\sum_{n=1}^\infty E[\Vert\phi(x_n)1_{[x_n\in U]} - \phi(b_n)\Vert^2]\ <\ \infty$, where $\Vert\cdot\Vert$ is
the Euclidean norm on $\bR^d$.
\end{theo}

  Note that under (G1), (G3) is equivalent to $\sum_{n=1}^\infty E[\Vert\phi(x_n)-\phi(b_n)\Vert^2
1_{[x_n\in U]}]<\infty$.

  The proof of Theorem~\ref{th1} will begin in the next section. Note that by Kolmogorov's 0\,-1 law, the independent
product $\hat{x}_n$ either converges almost surely or diverges almost surely.

When $G=\bR^d$ as an additive group, one may take $\phi$ to be the identity map on $\bR^d$ and $U$ to be the ball of
radius $r>0$ centered at $0$, then Theorem~\ref{th1} becomes precisely Kolmogorov's three-series theorem on $\bR^d$.

We briefly comment on the relation between the almost sure convergence and the convergence in distribution.
On Euclidean spaces, it is well known that the two convergences are equivalent for a series of independent random
variables. This is not true for an independent product on a Lie group $G$. Because if $G$ has
a compact subgroup $H\neq\{e\}$, then for any sequence of independent random variables $x_n$, each is distributed
according to the normalized Haar measure on $H$, the product $x_1x_2\cdots x_n$ converge in distribution to $x_1$, but
it is clearly not convergent almost surely. By Theorem~2.2.16\,(ii) in Heyer \cite{heyer}, if the only compact subgroup
of $G$ is $\{e\}$, then the convergence in distribution
and the almost sure convergence are equivalent for an infinite product of independent random variables in $G$.

For $k\geq 1$, let $\myM_k$ be the space of $k\times k$ real matrices, which may be identified
with $\bR^d$, where $d=k^2$. The Euclidean norm of $x=\{x_{ij}\}\in\myM_k$ is $\Vert x\Vert=\sqrt{\sum_{i,j}x_{ij}^2}$,
and it satisfies $\Vert xy\Vert\leq\Vert x\Vert\Vert y\Vert$ for $x,y\in\myM_k$.

Let $G$ be the group of $k\times k$ real matrices of nonzero determinants under matrix product. Its identity element $e$ is
the identity matrix $I$. Its Lie algebra is $\myM_k$ with the Lie group exponential map $\exp(x)$ being the usual matrix exponential $e^x=I+\sum_{n=1}^\infty
x^n/n!$.

\begin{theo} \label{th2}
Let $G$ be the matrix group as above, and let $x_n$ be a sequence of independent random variables in $G$.
Fix $r\in (0,\,1)$. Then $\hat{x}_n=x_1x_2\cdots x_n$ converges almost surely to a random matrix in $G$ if and only if the following three conditions hold.
\vspace{2ex}

\noindent (M1) \ $\sum_{n=1}^\infty P(\Vert x_n-I\Vert>r)\ <\ \infty$;

\noindent (M2) \ $b_1b_2\cdots b_n$ converges in $G$ as $n\to\infty$, where $b_n=I+E[(x_n-I)1_{[\Vert x_n-I\Vert\leq r]}]$; and

\noindent (M3) \ $\sum_{n=1}^\infty E(\Vert x_n-b_n\Vert^21_{[\Vert x_n-I\Vert\leq r]})\ <\ \infty$.

\end{theo}

\noindent {\bf Proof:} \ For $x\in G$, let $U=\{x\in G$; $\Vert x-I\Vert\leq r\}$ and $\phi(x)=x-I\in\myM_k$.
If $\Vert y\Vert<1$, then $I+y$ is invertible with $(I+y)^{-1}=I+\sum_{p=1}^\infty(-1)^py^p$. 
It follows that $\phi$ maps $U$ diffeomorphically onto
the ball of radius $r$ centered at $0$ in $\myM_k\equiv\bR^d$, and hence $\phi$ and $U$ satisfy the required properties.
Theorem~\ref{th1} may be applied with $b_n$ in (M2) being the $U$-truncated mean of $x_n$.
(G1) and (G2) are just (M1) and (M2), and (G3) is $\sum_nE[\Vert (x_n-I)1_{H_n}-(b_n-I)\Vert^2]<\infty$,
where $H_n=[\Vert x_n-I\Vert\leq r]$. Because $E[\Vert (x_n-I)1_{H_n}-(b_n-I)\Vert^2]=E[\Vert x_n-b_n\Vert^2 1_{H_n}]+\Vert b_n-I\Vert^2P(H_n^c)$,
by (M1), (G3) is equivalent to (M3). \ $\Box$
\vspace{2ex}

\noindent {\bf Example 1:} \ Let $y_n$ be a sequence of independent random variables in $\myM_k\equiv\bR^d$, $d=k^2$.
Assume $x_n=I+y_n$ is almost surely invertible. Note that this holds if $y_n$ has a continuous distribution. Also assume that for some $r\in (0,\,1)$, $E(y_n1_{[\Vert y_n\Vert\leq r]})=0$ for all $n$. Then $\hat{x}_n=x_1x_2\cdots x_n$ converges to an invertible random matrix $x_\infty$ almost surely if
\begin{equation}
  \sum_{n=1}^\infty E(\Vert y_n\Vert^2) <\ \infty. \label{sumEyn2}
\end{equation}
   To prove this claim, note that $b_n$ in (M2) is $I$ and (M2) holds trivially. Now (M1) is $\sum_{n=1}^\infty
P(\Vert y_n\Vert>r)<\infty$ and (M3) is $\sum_{n=1}^\infty E[\Vert y_n\Vert^21_{[\Vert y_n\Vert\leq r]}]<\infty$.
Because $P(\Vert y_n\Vert>r)\leq E(\Vert y_n\Vert^2)/r^2$, so (M1) and (M3) are implied by (\ref{sumEyn2}).
By Theorem~\ref{th2}, $\hat{x}_n$ converges almost surely in the matrix group $G$.
\vspace{2ex}

\noindent {\bf Example 2:} \ Let $y_n$ be independent random variables in $\myM_k\equiv\bR^d$, $d=k^2$. Assume $y_n$ is normal
of mean $0$. Then $\hat{x}_n=(I+y_1)\cdots (I+y_n)$ converges almost surely in the matrix group $G$ if and only if
(\ref{sumEyn2}) holds. To prove this, note that by the symmetry of a normal distribution, $E(y_n1_{[\Vert y_n\Vert\leq r]})=0$ for all $r>0$. By Example~1, (\ref{sumEyn2}) is a sufficient condition for the almost sure convergence of $\hat{x}_n$
in $G$. To see it is also necessary, it suffices to show that (\ref{sumEyn2}) is implied
by $\sum_n E[\Vert y_n\Vert^21_{[\Vert y_n\Vert\leq r]}]<\infty$ and $\sum_nP(\Vert y_n\Vert>r)<\infty$. This can
be done by an elementary computation of the normal distribution.
\vspace{2ex}

\noindent {\bf Example 3:} \ Let $y_n$ be a sequence of independent random variables in $\myM_k\equiv\bR^d$, $d=k^2$.
Assume there is $r>0$, which may be chosen arbitrarily small, such that $E(y_n1_{[\Vert y_n\Vert\leq r]})=0$
for all $n$. Then $\exp(y_1)\exp(y_2)\cdots\exp(y_n)$ converges in the matrix group $G$ almost surely if (\ref{sumEyn2}) holds. To prove this, apply Theorem~\ref{th1} to $x_n=\exp(y_n)$ with $\phi=\exp^{-1}$ on $U$, where $U$ is the diffeomorphic image of a small ball in $\myM_k\equiv\bR^d$
under $\exp$. The conditions may be verified as in Example~1.

\section{Sufficiency} \label{sec2}

For any sequence of independent random variables $x_n$ in $G$, by the Borel-Cantelli Lemma, if (G1) holds, then
almost surely, $x_n\in U$ except for finitely many $n$. On the other hand, if $\hat{x}_n=
x_1x_2\cdots x_n$ converges almost surely, then because $x_n=\hat{x}_{n-1}^{-1}\hat{x}_n\to e$, (G1) follows
from the Borel-Cantelli Lemma. Set $x_n'=x_n$ on $[x_n\in U]$ and $x_n'=e$ on $[x_n\in U^c]$. Then the
almost sure convergence of $x_1x_2\cdots x_n$ is equivalent to that of $x_1'x_2'\cdots x_n'$ and (G1). Note that $\phi(x_n)
1_{[x_n\in U]}=\phi(x'_n)=\phi(x'_n)1_{[x'_n\in U]}$, and all quantities in (G2) and (G3) (including $b_n$)
only depend on the restriction of $x_n$ on $U$. Therefore, (G2) and (G3) hold for $x_n$ if and only if they hold
for $x_n'$. Thus, as noted in \cite{maksimov}, to prove Theorem~\ref{th1}, we may, and will, assume all $x_n\in U$, and
prove that $\hat{x}_n$ converges almost surely in $G$ if and only if (G2) and (G3) hold.

We will prove the sufficiency part of Theorem~\ref{th1} in this section, and so assume (G2) and (G3).
Let $\mu_n$ be the distribution of $x_n$. Because $x_n\in U$, the $U$-truncated mean $b_n$ of $x_n$ is defined
by $\phi(b_n)=\mu_n(\phi)$, where $\mu_n(\phi)=\int\phi d\mu_n=E[\phi(x_n)]$. Set $\hat{x}_0=\hat{b}_0=e$.
For $n\geq 1$, let $z_n=\hat{b}_{n-1}x_nb_n^{-1}\hat{b}_{n-1}^{-1}$ and $\hat{z}_n=z_1z_2\cdots z_n$, and set $\hat{z}_0=e$.
It is easy to show by a simple induction on $n$ that for all $n\geq 0$,
\begin{equation}
\hat{x}_n = \hat{z}_n\hat{b}_n. \label{hatxhatzhatb}
\end{equation}
By (G2), it suffices to show that $\hat{z}_n$ converges in $G$ almost surely.

  Note that for $G=\bR^d$, $z_n$ is just the centered term $x_n-b_n$, and $\hat{z}_n=\hat{x}_n-\hat{b}_n$ is the sum
of the centered terms. To have $\hat{x}_n=\hat{z}_n\hat{b}_n$ on a non-commutative multiplicative Lie group $G$, $z_n$
has to be defined in the above rather complicated form.

By the lemma below, the almost sure convergence of $\hat{z}_n$ is equivalent to $z_mz_{m+1}\cdots z_n\to e$ almost surely
as $m\to\infty$ with $m<n$.

\begin{lemm} \label{lehatumn}
Let $u_n$ be independent random variables in $G$. Then $u_1u_2\cdots u_n$ converges almost surely as $n\to\infty$ if and only
if $u_mu_{m+1}\cdots u_n\to e$ almost surely as $m\to\infty$ with $m<n$.
\end{lemm}

\noindent {\bf Proof:} \ This is an easy consequence of the existence of a complete metric on $G$ that is invariant under
left translations and is compatible with the topology on $G$. The metric can be any left invariant Riemannian metric
on $G$. \ $\Box$
\vspace{2ex}

For any $f\in C_c^\infty(G)$, let $M_0f=f(e)$ and for $n\geq 1$, let
\begin{equation}
M_nf = f(\hat{z}_n) - \sum_{p=1}^n\int[f(\hat{z}_{p-1}\hat{b}_{p-1}xb_p^{-1}\hat{b}_{p-1}^{-1}) - f(\hat{z}_{p-1})]\mu_p(dx). \label{Mnf}
\end{equation}

\begin{lemm} \label{leMnf}
Let $\myF_n$ be the $\sigma$-algebra generated by $x_1,x_2,\ldots,x_n$. Then $E[M_nf\mid\myF_m]=M_mf$ for $m<n$, that
is, $M_nf$ is a martingale under the filtration $\{\myF_n\}$.
\end{lemm}

\noindent {\bf Proof:} \ Because $x_n$ are independent, for $m<p$,
\begin{eqnarray*}
&& E[\int f(\hat{z}_{p-1}\hat{b}_{p-1}xb_p^{-1}\hat{b}_{p-1}^{-1})\mu_p(dx)\mid\myF_m] \ = \ E[\int f(\hat{z}_mz_{m+1}\cdots z_{p-1}\hat{b}_{p-1}xb_p^{-1}\hat{b}_{p-1}^{-1})\mu_p(dx)\mid\myF_m] \\
&& \ \ = \ E[f(\hat{z}z_{m+1}\cdots z_{p-1}z_p)]\mid_{\hat{z}=\hat{z}_m} \ = \ E[f(\hat{z}_p)\mid\myF_m].
\end{eqnarray*}
Then $E[\int[f(\hat{z}_{p-1}\hat{b}_{p-1}xb_p^{-1}\hat{b}_{p-1}^{-1})-f(\hat{z}_{p-1})]\mu_p(dx)\mid\myF_m]=0$,
and $E[M_nf\mid\myF_m]=M_mf$. \ $\Box$
\vspace{2ex}

Fix an integer $m>0$ and a neighborhood $V$ of $e$. Let $f\in C_c^\infty(G)$ be such that $0\leq f\leq 1$, $f(e)=1$ and $f(x)=0$ for $x\in V^c$.
For $g\in G$, let $l_g$ be the left translation $x\mapsto gx$ on $G$, and let $f_m=f\circ l_{\hat{z}_m^{-1}}$. Let $\Lambda(m,V)$ be the event
that there is $n>m$ such that $z_{m+1}z_{m+2}\cdots z_n\in V^c$. To estimate $P[\Lambda(m,V)]$, let $\tau$ be the first time $n>m$ such
that $z_{m+1}z_{m+2}\cdots z_n\in V^c$ and set $\tau=\infty$ if $z_{m+1}z_{m+2}\cdots z_n\in V$ for all $n>m$. Then
\begin{equation}
P[\Lambda(m,V)] = E\{[f_m(\hat{z}_m)-f_m(\hat{z}_\tau)]1_{\Lambda(m,V)}\} = \lim_{n\to\infty}E\{[f_m(\hat{z}_m)-f_m(\hat{z}_{\tau\wedge n})]
1_{\Lambda(m,V)}\}, \label{PLambdam}
\end{equation}
where $\tau\wedge n=\min(\tau,n)$. Because $E\{[f_m(\hat{z}_m)-f_m(\hat{z}_{\tau\wedge n})]1_{\Lambda(m,V)}\}\leq E[1-f_m(\hat{z}_{\tau\wedge n})]=
E[f_m(\hat{z}_m)-f_m(\hat{z}_{\tau\wedge n})]$ and $E[M_{\tau\wedge n}f_m]=E\{E[M_{\tau\wedge n}f_m\mid\myF_m]\}=E[M_mf_m]$,
\begin{eqnarray}
E\{[f_m(\hat{z}_m)-f_m(\hat{z}_{\tau\wedge n})]1_{\Lambda(m,V)}\} &\leq& -E\{\sum_{p=m+1}^{\tau\wedge n}
\int[f(\hat{z}_{p-1}\hat{b}_{p-1}xb_p^{-1}\hat{b}_{p-1}^{-1}) - f(\hat{z}_{p-1})]\mu_p(dx)\} \nonumber \\
\leq && \hspace{-0.3in} \sum_{p=m}^\infty E\{\vert\int[f(\hat{z}_{p-1}\hat{b}_{p-1}xb_p^{-1}\hat{b}_{p-1}^{-1}) - f(\hat{z}_{p-1})]\mu_p(dx)\vert\} \label{sum}.
\end{eqnarray}

We will write $\hat{z},\hat{b},b,\mu$ for $\hat{z}_{p-1},\hat{b}_{p-1},b_p,\mu_p$ for simplicity. For $x\in U$, by the Taylor expansion
of $f(\hat{z}\hat{b}xb^{-1}\hat{b}^{-1})=f(\hat{z}\hat{b}\phi^{-1}(\phi(x))b^{-1}\hat{b}^{-1})$ at $x=b$, noting $\mu(U^c)=0$,
\begin{equation}
\int [f(\hat{z}\hat{b}xb^{-1}\hat{b}^{-1}) - f(\hat{z})]\mu(dx) = \int\{\sum_i f_i(\hat{z},\hat{b},b)[\phi_i(x)-\phi_i(b)]\}\mu(dx) + r,
\label{inttaylor}
\end{equation}
where
\begin{equation}
f_i(\hat{z},\hat{b},b)\ =\ \frac{\partial}{\partial\phi_i}
f(\hat{z}\hat{b}\phi^{-1}(\phi(x))b^{-1}\hat{b}^{-1})\mid_{x=b} \label{fi}
\end{equation}
and the remainder $r$ satisfies $\vert r\vert\leq c\mu(\Vert\phi-\phi(b)\Vert^2)$
for some constant $c>0$. Because $\mu(\phi_i)=\phi_i(b)$, $\int[\phi_i(x)-\phi_i(b)]\mu(dx)=0$, and then by (\ref{inttaylor}),
\begin{equation}
\vert\int [f(\hat{z}\hat{b}xb^{-1}\hat{b}^{-1}) - f(\hat{z})]\mu(dx)\vert\ =\ \vert r\vert\ \leq\ c\mu(\Vert\phi-\phi(b)\Vert^2). \label{intfleq}
\end{equation}

It now follows from (\ref{PLambdam}) and (\ref{sum}) that $P[\Lambda(m,V)]\leq c\sum_{n=m}^\infty\mu_n(\Vert\phi-\phi(b_n)\Vert^2)$.
Let $\varepsilon\in (0,\,1)$ and let $V_k$ be a sequence of neighborhoods of $e$ with $V_k\downarrow\{e\}$ as $k\uparrow\infty$. By (G3),
for each $k\geq 1$, there is an integer $m_k$ such that $P[\Lambda(m_k,V_k)]<\varepsilon^k$. Then $\sum_{k=1}^\infty P[\Lambda(m_k,V_k)]\leq\sum_{k=1}^\infty\varepsilon^k=\varepsilon/(1-\varepsilon)$.
By Lemma~\ref{lehatumn}, $P(\hat{z}_n$ converges$)\geq P[\cap_{k=1}^\infty\Lambda(m_k,V_k)^c]\geq 1-\sum_{k=1}^\infty P[\Lambda(m_k,V_k)]\geq 1-\varepsilon/(1-\varepsilon)\to 1$ as $\varepsilon\to 0$. This proves $\hat{z}_n$ converges
almost surely.

\section{Necessity, part 1} \label{sec3}

We will now prove (G2) and (G3) under the assumption that $\hat{x}_n$ converges almost surely and all $x_n\in U$. This proof
is more complicated and will require another section.

Because $x_n=\hat{x}_{n-1}^{-1}\hat{x}_n\to e$ almost surely, by the Borel-Cantelli Lemma,
\begin{equation}
\forall\ {\rm neighborhood}\ V\ {\rm of}\ e,\ \ \ \ \sum_{n=1}^\infty P(x_n\in V^c)\ <\ \infty. \label{strongerthana}
\end{equation}

We also have
\begin{equation}
b_n\to e\ \ \ \ {\rm and} \ \ \ \ \mu_n(\Vert\phi-\phi(b_n)\Vert^2)\to 0\ \ \ \ {\rm as}\ \ n\to\infty. \label{bntoe}
\end{equation}

For $m<n$, let $\hat{x}_{m,n}=x_{m+1}x_{m+2}\cdots x_n$ and $\hat{b}_{m,n}=b_{m+1}b_{m+2}\cdots b_n$. If either (G2)
or (G3) does not hold, then there are a neighborhood $V$ of $e$, $\varepsilon>0$ and two sequences of integers $m_k$
and $n_k$ with $VV\subset U$, $m_k<n_k$ and $m_k\uparrow\infty$ as $k\uparrow\infty$ such that for each $k\geq 1$,
\[\mbox{either \ $\sum_{p=m_k+1}^{n_k}\mu_p(\Vert\phi-\phi(b_p)\Vert^2)\geq\varepsilon$ \ or \ $\hat{b}_{m_k,\,n_k}\in V^c$.}\]
Because of (\ref{bntoe}), by choosing $m_1$ large enough, we have $b_n\in V$ and $\mu_n(\Vert\phi-\phi(b_n)\Vert^2)\leq
\varepsilon$ for $n>m_1$. Thus, by suitably reducing $n_k$, we obtain that for each $k\geq 1$, either
\vspace{2ex}

\noindent (i) \ $\varepsilon\leq\sum_{p=m_k+1}^{n_k}\mu_p(\Vert\phi-\phi(b_p)\Vert^2)\leq
2\varepsilon$, and $\hat{b}_{m_k,\,p}\in U$ for $m_k<p\leq n_k$; or

\noindent (ii) \ $\sum_{p=m_k+1}^{n_k}\mu_p(\Vert\phi-\phi(b_p)\Vert^2)\leq
2\varepsilon$, $\hat{b}_{m_k,\,n_k}\in V^c$, and $\hat{b}_{m_k,\,p}\in U$ for $m_k<p\leq n_k$.
\vspace{2ex}

We will derive a contradiction from either (i) or (ii) above. We will embed the
partial products $x_{m_k,\,p}$ and $b_{m_k,\,p}$, for $m_k<p\leq n_k$, into a process $\tilde{x}_t^k$ and
a function $\tilde{b}_t^k$ on $[0,\,1]$ respectively. The main idea is to obtain a martingale property
for the process $\tilde{z}_t^k$, defined by $\tilde{x}_t^k=\tilde{z}_t^k\tilde{b}_t^k$, similar to the martingale
property for $\hat{z}_n$ in the last section, to show the limit $\tilde{z}_t$ of $\tilde{z}_t^k$ satisfies
an integral equation, and then to derive a contradiction. This is similar to the approaches
in \cite{feinsilver,Liao} for processes in Lie groups with independent increments.

Let $\gamma_k$ be a strictly increasing function from $\{m_k,m_k+1,\ldots,n_k\}$ into $[0,\,1]$ with $\gamma_k(m_k)=0$
and $\gamma_k(n_k)=1$. Let $t_{k,p}=\gamma_k(p)$ for $m_k\leq p\leq n_k$. Then $t_{k,m_k}=0$ and $t_{k,n_k}=1$.
Let $\tilde{x}_t^k=\tilde{b}_t^k=e$ for $0\leq t<t_{k,\,m_k+1}$. For $m_k<p<n_k$ and $t_{k,p}\leq t<t_{k,\,p+1}$, let
\begin{equation}
\tilde{x}_t^k = \hat{x}_{m_k,\,p}\ \ \ \ {\rm and}\ \ \ \ \tilde{b}_t^k = \hat{b}_{m_k,\,p}. \label{tildexb}
\end{equation}
  Set $\tilde{x}_t^k=\hat{x}_{m_k,\,n_k}$ and $\tilde{b}_t^k=\hat{b}_{m_k,\,n_k}$ for $t\geq 1$.
Then $\tilde{x}_t^k$ and $\tilde{b}_t^k$ are respectively a step process and a step function, which are
right continuous with jumps $x_p$ and $b_p$ at $t=t_{k,p}$.

  Note that by Lemma~\ref{lehatumn}, almost surely, $\tilde{x}_t^k\to e$ as $k\to\infty$ uniformly in $t$.

A continuous function $A(t)=\{A_{ij}(t)\}$ from $\bR_+=[0,\,\infty)$ to the space of $d\times d$ symmetric real matrices
is called a covariance matrix function if $A(0)=0$ and for $s<t$, $A(t)-A(s)\geq 0$ (nonnegative definite). Let
\begin{equation}
A_{ij}^k(t) = \sum_{0<t_{k,p}\leq t}\int_G[\phi_i(x)-\phi_i(b_p)][\phi_j(x)-\phi_j(b_p)]\mu_p(dx). \label{Aijkt}
\end{equation}
Then $A^k(t)=\{A_{ij}^k(t)\}$ is almost a covariance matrix function except that it is not continuous, but $A^k(t)=A^k(1)$
for $t\geq 1$. Let $Q^k(t)$ be the trace of $A^k(t)$. Then
\begin{equation}
Q^k(t) = \sum_{0<t_{k,p}\leq t}\mu_p(\Vert\phi-\phi(b_p)\Vert^2), \label{Qkt}
\end{equation}
and for $s<t$,
\begin{equation}
\vert A_{ij}^k(t)-A_{ij}^k(s)\vert\ \leq\ Q^k(t)-Q^k(s). \label{AleqQ}
\end{equation}
  Note that $Q^k(t)$ is a nondecreasing step function in $t$ with a jump $\mu_p
(\Vert\phi-\phi(b_p)\Vert^2)$ at $t=t_{k,p}$, $Q^k(t)=0$
for $0\leq t< t_{m_k,\,m_k+1}$ and $Q^k(t)=Q^k(1)=\sum_{m_k<p\leq n_k}\mu_p(\Vert\phi-\phi(b_p)\Vert^2)$ for $t\geq 1$.
By either (i) or (ii), $Q^k(t)\leq 2\varepsilon$, and by (\ref{bntoe}), the jumps of $Q^k(t)$ converge to $0$ uniformly in $t$ as $k\to\infty$. It follows that the function $\gamma_k$ may be chosen properly such that
\begin{equation}
Q^k(t) - Q^k(s)\ \leq\ 2\varepsilon (t-s) + \varepsilon_k,\ \ \ \ 0\leq s<t\leq 1, \label{Qktequicont}
\end{equation}
where $\varepsilon_k\to 0$ as $k\to\infty$. Roughly speaking, this means the functions $Q^k(t)$ are equi-continuous
for large $k$. Because of (\ref{bntoe}), by either (i) or (ii), $n_k-m_k\to\infty$ as $k\to\infty$, and hence $\gamma_k$ may
be chosen to satisfy, besides (\ref{Qktequicont}),
\begin{equation}
\max_{p>m_k+1}(t_{k,p}-t_{k,\,p-1})\ \to\ 0\ \ \ \ {\rm as}\ \ k\to\infty. \label{tt0}
\end{equation}

\begin{lemm} \label{le1}
There is a covariance matrix function $A(t)$ with $A(t)=A(1)$ for $t\geq 1$ such that along a subsequence
of $k\to\infty$, $A^k(t)\to A(t)$ for any $t\geq 0$.
\end{lemm}

\noindent {\bf Proof} \ Let $\Lambda$ be a countable dense subset of $[0,\,1]$. Under either (i)
or (ii), $Q^k(t)$ is bounded. By (\ref{AleqQ}), along a subsequence of $k\to\infty$, $A^k(t)$ converges for any $t\in\Lambda$.
By (\ref{Qktequicont}), the convergence holds for all $t\geq 0$, and $A(t)$ is continuous in $t$. \ $\Box$
\vspace{2ex}

Let $Y$ be a smooth manifold equipped with a compatible metric $\rho$ and let $y$: $[0,\,1]\to Y$ be a continuous function.
For each $k$, let $y^k$: $[0,\,1]\to Y$ be a step function that is constant on $[t_{k,\,p-1},\,t_{k,p})$ for
each $p=m_k+1,\ldots,n_k$. Assume for any $t>0$, $\rho(y^k(t_{k,p}),y(t_{k,p}))\to 0$ as $k\to\infty$ uniformly
for $t_{k,p}\leq t$. Let $F(y,g)=\{F_{ij}(y,g)\}$ be a bounded continuous matrix-valued function on $Y\times G$.

\begin{lemm} \label{le2}
Assume the above and let $A(t)$ be the covariance matrix function in Lemma~\ref{le1}. Then for any $t>0$, along the subsequence of $k\to\infty$ in Lemma~\ref{le1},
\begin{eqnarray}
&&  \sum_{0<t_{k,p}\leq t}\,\sum_{i,j=1}^d\int_G F_{ij}(y^k(t_{k,\,p-1}),b_p)[\phi_i(x)-\phi_i(b_p)][\phi_j(x)-\phi_j(b_p)]\mu_p(dx) \nonumber \\
&\to& \sum_{i,j=1}^d \int_0^t F_{ij}(y(s),e) dA_{ij}(s). \label{fjkyntni1to}
\end{eqnarray}
\end{lemm}

\noindent {\bf Proof} \ By the uniform convergence $\rho(y^k(t_{k,p}),y(t_{k,p}))\to 0$, $F(y^k(t_{k,p}),b)-
F(y(t_{k,p}),b)\to 0$ as $k\to\infty$ uniformly for $t_{k,p}\leq t$ and for $b$ in a compact set.
Because when $k\to\infty$, $b_p\to e$ uniformly for $p>m_k$,
we may replace $y^k$ and $b_p$ by $y$ and $e$ in the proof.

Let $r>0$ be an integer. For any two expressions $A$ and $B$ depending on $(k,r)$,
we will write $A\approx B$ if $\vert A-B\vert\to 0$ as $r\to\infty$ uniformly in $k$. Then
\begin{eqnarray*}
&& \sum_{0<t_{k,p}\leq t}\,\sum_{i,j=1}^d\int_G F_{ij}(y(t_{k,\,p-1}),e)
[\phi_i(x)-\phi_i(b_p)][\phi_j(x)-\phi_j(b_p)]\mu_p(dx) \\
&\approx& \sum_{i,j=1}^d\sum_{q=0}^{r-1}\,\sum_{qt/r<t_{k,p}\leq (q+1)t/r}\int_G F_{ij}(y(\frac{qt}{r}),e)
[\phi_i(x)-\phi_i(b_p)][\phi_j(x)-\phi_j(b_p)]\mu_p(dx) \\
&& \mbox{(where $\sum_{qt/r<t_{k,p}\leq (q+1)t/r}(\cdots)=0$ if $(\frac{qt}{r},\,\frac{(q+1)t}{t}]$ contains no $t_{k,p}$)} \\
&\to& \sum_{i,j=1}^d\sum_{q=0}^{r-1} F_{ij}(y(\frac{qt}{r}),e)[A_{ij}(\frac{(q+1)t}{r})-A_{ij}(\frac{qt}{r})]
\ \ \ \ \mbox{(as $k\to\infty$, by Lemma~\ref{le1})} \\
&\approx& \sum_{i,j=1}^d\int_0^t F_{ij}(y(s),e)dA_{ij}(s). \ \ \ \ \Box
\end{eqnarray*}

We now define a new process $\tilde{z}_t^k$, similar in the way as the sequence $z_n$ is defined from $x_n$ and $b_n$
in \S\ref{sec2}, by setting $\tilde{z}_t^k=e$ for $0\leq t<t_{k,\,m_k+1}$, and inductively
\begin{equation}
\tilde{z}_t^k\ =\ \tilde{z}_{t_{k,\,p-1}}^k\tilde{b}_{t_{k,\,p-1}}^kx_pb_p^{-1}(\tilde{b}_{t_{k,p-1}}^k)^{-1} \label{tildeztk}
\end{equation}
for $t_{k,p}\leq t<t_{k,\,p+1}$, $p=m_k+1,\ldots,n_k$, setting $t_{k,\,n_k+1}=\infty$ here. Then $\tilde{z}_t=\tilde{z}_1$
for $t>1$, and a simple induction on $p$ shows that $\tilde{x}_t^k=\tilde{z}_t^k\tilde{b}_t^k$ for all $t\geq 0$.

For $f\in C_c^\infty(G)$, let $\tilde{M}_t^kf=f(\tilde{z}_t^k)=f(e)$ for $0\leq t<t_{k,\,m_k+1}$, and let
\begin{equation}
\tilde{M}_t^kf\ =\ f(\tilde{z}_t^k) - \sum_{0<t_{k,p}\leq t}\int_G[f(\tilde{z}_{t_{k,\,p-1}}^k\tilde{b}_{t_{k,\,p-1}}^k xb_p^{-1}(\tilde{b}_{t_{k,\,p-1}}^k)^{-1}) - f(\tilde{z}_{t_{k,\,p-1}})]\mu_p(dx), \label{tildeMtkf}
\end{equation}
for $t\geq t_{k,\,m_k+1}$.

\begin{lemm} \label{letildeMtkf} \ $\tilde{M}_t^kf$ is a martingale under the natural filtration of process $\tilde{z}_t^k$.
\end{lemm}

\noindent {\bf Proof:} \ This is proved in the same way as in Lemma~\ref{leMnf} for $M_nf$ to be a martingale. \ $\Box$
\vspace{2ex}

Because $\tilde{x}_t^k=\tilde{z}_t^k\tilde{b}_t^k$ and $\tilde{x}_t^k\to e$ uniformly in $t$ as $k\to\infty$
almost surely, if $\tilde{b}_t^k$ converges to some continuous path $\tilde{b}_t$ in $G$ uniformly in $t$
as $k\to\infty$, then $\tilde{z}_t^k\to\tilde{z}_t=\tilde{b}_t^{-1}$ uniformly in $t$ almost surely. This will be assumed
in the rest of this section.

By a computation using Taylor expansion similar to the one in the last section, but up to the second order, noting the integrals of the first order terms vanish as before,
\[\tilde{M}_t^kf\ =\ f(\tilde{z}_t^k) - \sum_{0<t_{k,p}\leq t} \sum_{i,j}\int_G f_{ij}(\tilde{z}_{t_{k,\,p-1}}^k,
\tilde{b}_{t_{k,\,p-1}}^k,b_p)[\phi_i(x)-\phi_i(b_p)][\phi_j(x)-\phi_p(b_p)]\mu_p(dx) + r_k,\]
where
\[f_{ij}(\tilde{z},\tilde{b},b)\ =\ \frac{\partial^2}{\partial\phi_i\partial\phi_j}
f(\tilde{z}\tilde{b}\phi^{-1}(\phi(x))b^{-1}\tilde{b}^{-1})\mid_{x=b},\]
and the reminder $r_k$ may be divided into an integral over a small neighborhood $V$ of $e$ and an integral
over $V^c$. The former is controlled by $c_VQ^k(t)\leq c_V(2\varepsilon)$, where the constant $c_V\to 0$ as $V\downarrow\{e\}$, and the latter is controlled by $\sum_{m_k<p\leq n_k}\mu_p(V^c)$ which converges to $0$ as $k\to\infty$ by (\ref{strongerthana}). Therefore, $r_k\to 0$ as $k\to\infty$. By Lemma~\ref{le2} with $Y=G\times G$
and $y^k(t)=(\tilde{z}_t^k,\tilde{b}_t^k)\to y(t)=(\tilde{z}_t, \tilde{b}_t)$,
it follows that $\tilde{M}_t^kf$ converges to the martingale
\[\tilde{M}_tf = f(\tilde{z}_t)-\sum_{i,j}\int_0^tf_{ij}(\tilde{z}_s,\tilde{b}_s,e)dA_{ij}(s)\]
as $k\to\infty$. Because $\tilde{z}_t=\tilde{b}_t^{-1}$ is non-random, the martingale $\tilde{M}_tf$ must be $f(e)$, and then
for any $f\in C_c^\infty(G)$ with $f(e)=0$,
\begin{equation}
f(\tilde{z}_t) \ = \ \sum_{i,j}\int_0^t[\frac{\partial^2}{\partial\phi_i\,\partial\phi_j}f(\phi^{-1}(\phi(x))\tilde{z}_s)\mid_{x=e}]dA_{ij}(s). \label{ftildezt}
\end{equation}

Let $t_0$ be the largest nonnegative real number $\leq 1$ such that $\tilde{z}_s=e$ and $A(s)=0$ for $s\leq t_0$.
We will show $t_0=1$. Suppose $t_0<1$. Then (\ref{ftildezt}) holds for $t\geq t_0$ with $\int_0^t$ replaced by $\int_{t_0}^t$.
Without loss of generality, we will assume $t_0=0$. Substitute $f=\phi_\beta^2$ in (\ref{ftildezt}), then the integrand
is $2\delta_{i\beta}\delta_{j\beta}+\varepsilon_s$, where $\varepsilon_s$ denotes any function satisfying $\varepsilon_s\to 0$ as $s\to 0$. It follows that $\phi_\beta(\tilde{z}_t)^2=2A_{\beta\beta}(t)+\varepsilon_tT_t$,
where $T_t=\trace[A(t)]$. Then $\Vert\phi(\tilde{z}_t)\Vert^2=2T_t+\varepsilon_tT_t$. Now
let $f=\phi_\beta$ and then (\ref{ftildezt}) yields $\phi_\beta(\tilde{z}_t)=\varepsilon_tT_t$. This implies $\vert
\phi_\beta(\tilde{z}_t)\vert\leq c\Vert\phi(\tilde{z}_t)\Vert^2$ for some constant $c>0$, which is clearly
impossible. This shows that $t_0=1$, and hence $\tilde{z}_t=e$ and $A(t)=0$ for all $t\geq 0$.

If (i) holds, then $\trace[A(1)]=\lim_kQ^k(1)=\lim_k\sum_{p=m_k+1}^{n_k}\mu_p(\Vert\phi-
\phi(b_p)\Vert^2)\geq\varepsilon$, which contradicts to $A(t)=0$. Thus (i) cannot hold. If (ii) holds,
then $\tilde{b}_1=\lim_k\tilde{b}_1^k=\lim_k\hat{b}_{m_k,n_k}$ belongs to the closure of $V^c$,
which contradicts to $\tilde{b}_t=\tilde{z}_t^{-1}=e$.
We have proved that neither (i) nor (ii) holds, and hence (G2) and (G3) must
hold, under the assumption that $\tilde{b}_t^k\to\tilde{b}_t$ as $k\to\infty$ uniformly in $t$
for some continuous path $\tilde{b}_t$ in $G$.

\section{Necessity, part 2}

It remains to show that $\tilde{b}_t^k\to\tilde{b}_t$ as $k\to\infty$ uniformly in $t$ for some continuous path $\tilde{b}_t$ in $G$. A rcll path is a right continuous path with left limits, and a process with rcll paths will be called a rcll process.
Let $D(G)$ be the space of rcll paths in $G$. Equipped with  the Skorohod metric, $D(G)$ is
a complete separable metric space (see \cite[chapter~3]{ek}).  A sequence of rcll processes $y_t^k$ in $G$
are said to converge weakly to a rcll process $y_t$ if $y_\cdot^k\to y_\cdot$ in distribution as $D(G)$-valued random variables. The sequence $y_t^k$ are called relatively weak compact in $D(G)$ if any subsequence has
a further subsequence that converge weakly.

We will show that $\tilde{z}_t^k$ are relatively weak compact. Let $V$ be a neighborhood of $e$. The amount of time it takes for a rcll process $y_t$ to make $V^c$-displacement from a stopping time $\sigma$ (under the natural filtration of
process $y_t$) is denoted as $\tau_V^\sigma$, that is,
\begin{equation}
\tau_V^\sigma = \inf\{t>0;\ \ y_\sigma^{-1}y_{\sigma+t}\in V^c\}\ \ \ \ \mbox{($\inf$ of an empty set is $\infty$).}
\label{tauUsigma}
\end{equation}
For a sequence of processes $y_t^k$ in $G$, let $\tau_V^{\sigma,k}$
be the $V^c$-displacement time for $y_t^k$ from $\sigma$.

The following lemma is Lemma~16 in \cite{Liao} and provides a criterion for the relative compactness. It is a slightly improved version of a lemma in \cite{feinsilver}.

\begin{lemm} \label{le3}
A sequence of rcll processes $y_t^k$ in $G$ are relatively weak compact in $D(G)$ if for any constant $T>0$ and
any neighborhood $V$ of $e$,
\begin{equation}
\overline{\lim_{k\to\infty}}\sup_{\sigma\leq T} P(\tau_V^{\sigma,k}
<\delta)\ \to\ 0\ \ {\rm as}\ \delta\to 0, \label{supPtauU}
\end{equation}
and
\begin{equation}
\overline{\lim_{k\to\infty}}\sup_{\sigma\leq T}
P[(y_{\sigma-}^k)^{-1} y_\sigma^k\in K^c]\ \to\ 0 \ \ \mbox{as
compact}\ K\uparrow G, \label{supPxsigma}
\end{equation}
where $\sup_{\sigma\leq T}$ is taken over all stopping times $\sigma\leq T$.
\end{lemm}

We will apply Lemma~\ref{le3} to $y_t^k=\tilde{z}_t^k$. Because $\tilde{z}_t^k=\tilde{z}_1^k$ for $t>1$, we may take $T=1$ in Lemma~\ref{le3}.
Let $f\in C_c^\infty(G)$ be such that $0\leq f\leq 1$ on $G$, $f(e)=1$ and $f=0$ on $V^c$. For any stopping time $\sigma\leq 1$,
write $\tau $ for $\tau_V^{\sigma,k}$ and let $f_\sigma=f\circ l_z$ with $z=(\tilde{z}_\sigma^k)^{-1}$. Then
\begin{equation}
P(\tau<\delta) = E[f_\sigma(\tilde{z}_\sigma^k)-f_\sigma(\tilde{z}_{\sigma+\tau}^k);\,\tau<\delta] \leq
E[f_\sigma(\tilde{z}_\sigma^k)-f_\sigma(\tilde{z}_{\sigma+\tau\wedge\delta}^k)], \label{Ptaudelta2}
\end{equation}
  noting $f_\sigma(z_\sigma^k)=1$, $f_\sigma(z_{\sigma+\tau}^k)=0$ and $\tau=\tau\wedge\delta$ on $[\tau<\delta]$.
Because $\tilde{M}_t^kf$ given by (\ref{tildeMtkf}) is a martingale for any $f\in C_c^\infty(G)$, and $\sigma$ and $\sigma+\tau\wedge\delta$
are stopping times,
\[E[\tilde{M}_\sigma^kf_\sigma - \tilde{M}_{\sigma+\tau\wedge\delta}^kf_\sigma] = E\{E[\tilde{M}_\sigma^kf_\sigma -
\tilde{M}_{\sigma+\tau\wedge\delta}^k f_\sigma\mid\myF_\sigma]\} = 0.\]
Writing $\tilde{z},\tilde{b},b,\mu$ for $\tilde{z}_{t_{k,\,p-1}}^k,\tilde{b}_{t_{k,\,p-1}}^k,b_p, \mu_p$, by (\ref{tildeMtkf}) and (\ref{Ptaudelta2}),
we obtain
\begin{eqnarray}
P(\tau<\delta) &\leq& - E\{\sum_{\sigma<t_{k,p}\leq\sigma+\tau\wedge\delta}\int_G
[f_\sigma(\tilde{z}\tilde{b}xb^{-1}\tilde{b}^{-1})-f_\sigma(\tilde{z})]\mu(dx)\} \nonumber \\
&\leq& E\{\sum_{\sigma<t_{k,p}\leq\sigma+\delta}\vert
\int_G[f_\sigma(\tilde{z}\tilde{b}xb^{-1}\tilde{b}^{-1})-f_\sigma(\tilde{z})]\mu(dx)\vert\}.
\label{Ptaudelta}
\end{eqnarray}
Performing the same computation leading to (\ref{intfleq}) shows that for some constant $c>0$,
\[P(\tau<\delta)\ \leq\ cE[Q^k(\sigma+\delta)-Q^k(\sigma)].\]
By (\ref{Qktequicont}), $E[Q^k(\sigma+\delta)-Q^k(\sigma)]\leq 2\varepsilon\delta+\varepsilon_k$.
It follows that $\overline{\lim}_{k\to\infty}\sup_{\sigma\leq 1}P(\tau<\delta)\leq 2c
\varepsilon\delta$. This shows that the condition (\ref{supPtauU}) is satisfied for $y_t^k=\tilde{z}_t^k$.

  To verify (\ref{supPxsigma}), note that because $\tilde{x}_t^k=\tilde{z}_t^k\tilde{b}_t^k$,
\[P[(\tilde{z}_{\sigma-}^k)^{-1}\tilde{z}_\sigma^k\in K^c] = P[(\tilde{x}_{\sigma-}^k)^{-1}\tilde{x}_\sigma^k \in
(\tilde{b}_{\sigma-}^k)^{-1}K^c \tilde{b}_\sigma^k].\]
By either (i) or (ii), $\tilde{b}_t^k$ are bounded in $k$, when $K$ is large, $(\tilde{b}_{\sigma-}^k)^{-1}K \tilde{b}_\sigma^k$ contains
a fixed neighborhood $H$ of $e$. Because $(\tilde{b}_{\sigma-}^k)^{-1}K^c \tilde{b}_\sigma^k=((\tilde{b}_{\sigma-}^k)^{-1}K\tilde{b}_\sigma^k)^c$,
it follows that
\[P[(\tilde{z}_{\sigma-}^k)^{-1}\tilde{z}_\sigma^k\in K^c]\ \leq\ P[(\tilde{x}_{\sigma-}^k)^{-1}\tilde{x}_\sigma^k \in H^c]
\ \leq\ \sum_{p>m_k}\mu_p(H^c) \to 0\]
  as $k\to\infty$. This verifies (\ref{supPxsigma}) even before taking $K\uparrow G$.

By Lemma~\ref{le3}, $\tilde{z}_t^k$ are relatively weak compact, and hence along a subsequence
of $k\to\infty$, $\tilde{z}_t^k$ converge weakly to a rcll process $\tilde{z}_t$ in $G$. As $D(G)$-valued random
variables, $\tilde{z}_\cdot^k$ converge in distribution to $\tilde{z}_\cdot$. It is well known (see
for example Theorem~1.8 in \cite[chapter~3]{ek}) that there are $D(G)$-valued
random variables $\tilde{z}_\cdot'^k$ and $\tilde{z}_\cdot'$, possibly on a different
probability space, such that $\tilde{z}_\cdot'$ is equal to $\tilde{z}_\cdot$ in distribution, $\tilde{z}_\cdot'^k$ is
equal to $\tilde{z}_\cdot^k$ in distribution for each $k$, and $\tilde{z}_\cdot'^k\to\tilde{z}_\cdot'$ almost surely.
Because $\tilde{x}_\cdot^k=\tilde{z}_\cdot^k\tilde{b}_\cdot^k\to e$
almost surely, where $e$ is regarded as a constant path in $G$, $\tilde{x}_\cdot'^k=\tilde{z}_\cdot'^k\tilde{b}_\cdot^k\to e$ in distribution. As the limit $e$ is non-random, $\tilde{x}_\cdot'^k\to e$ in probability. Then along a further subsequence
of $k\to\infty$, $\tilde{x}_\cdot'^k\to e$ almost surely, and hence $\tilde{b}_\cdot^k=
(\tilde{z}_\cdot'^k)^{-1}\tilde{x}_\cdot'^k\to (\tilde{z}_\cdot')^{-1}$.

The convergence $\tilde{b}_t^k\to\tilde{b}_t=(\tilde{z}_t')^{-1}$ under the Skorohod metric means (see Proposition~5.3(c) in \cite[chapter~3]{ek})
that there are continuous strictly increasing functions $\lambda_k$: $\bR_+\to\bR_+$ such that as $k\to\infty$, $\lambda_k(t)-t
\to 0$ and $r(\tilde{b}_t^k,\tilde{b}_{\lambda_k(t)})\to 0$ uniformly for $0\leq t\leq 1$, where $r$ is a compatible metric on $G$.
If $\tilde{b}_t$ has a jump of size $r(\tilde{b}_{s-},\tilde{b}_s)>0$ at time $s$, then $\tilde{b}_t^k$ would
have a jump of size close to $r(\tilde{b}_{s-}^\gamma,\tilde{b}_s)$ at time $t=\lambda_k^{-1}(s)$,
which is impossible because the jumps of $\tilde{b}_t^k$ are uniformly small
when $k$ is large. It follows that $\tilde{b}_t$ is continuous in $t$ and hence $\tilde{b}_t^k\to\tilde{b}_t$ uniformly in $t$ as $k\to\infty$.
\vspace{2ex}

\noindent {\bf Acknowledgement:} \ The author wishes to thank David Applebaum and an anonymous referee for some very
helpful comments.


\begin{thebibliography}{99}

\bibitem{billingsley} Billingsley, P. ``Probability and measure", second edition, Wiley (1986).

\bibitem{ek} Ethier, S.N. and Kurtz, T.G., ``Markov processes, characterization and convergence", John Wiley (1986).

\bibitem{feinsilver} Feinsilver, P., ``Processes with independent increments on a Lie group", Trans. Amer. Math. Soc. 242,
(1978), 73-121.

\bibitem{heyer} Heyer, H., ``Probability measures on locally compact groups", Springer (1977).


\bibitem{Liao} Liao, M., ``Inhomogeneous \levy processes in Lie groups and homogeneous spaces", J. Theor. Probab. 27, (2014),
315-357.

\bibitem{maksimov} Maksimov, V.M., ``The principle of convergence ``almost everywhere" in Lie groups", Math. USSR Sbornik 20,
(1973), 543-555.

\end{thebibliography}
\end{document}